\documentclass[leqno]{article}
\usepackage{amsmath,amssymb,amscd,latexsym}
\usepackage[all]{xy}
\newcommand{\C}{{\mathbb{C}}}
\newcommand{\F}{{\mathbb{F}}}
\newcommand{\oF}{\overline{\F}}
\newcommand{\Q}{{\mathbb{Q}}}
\newcommand{\R}{{\mathbb{R}}}
\newcommand{\Z}{{\mathbb{Z}}}
\newcommand{\abb}{\mathrm{ab}}
\newcommand{\car}{\mathrm{char}}
\newcommand{\ddet}{\mathrm{det}}
\newcommand{\ddis}{\mathrm{dis}}
\newcommand{\id}{\mathrm{id}}
\newcommand{\ord}{\mathrm{ord}}
\newcommand{\rank}{\mathrm{rank}}
\newcommand{\rk}{\mathrm{rk}\,}
\newcommand{\sgn}{\mathrm{sgn}\,}
\newcommand{\spec}{\mathrm{spec}\,}
\newcommand{\Imm}{\mathrm{Im}\,}
\newcommand{\Ker}{\mathrm{Ker}\,}
\newcommand{\RRe}{\mathrm{Re}\,}

\newcommand{\tors}{\mathrm{tors}}
\newcommand{\tr}{\mathrm{tr}}
\newcommand{\Ah}{{\mathcal A}}
\newcommand{\Dh}{{\mathcal D}}
\newcommand{\Fh}{{\mathcal F}}
\newcommand{\ea}{\mathfrak{a}}
\newcommand{\eb}{\mathfrak{b}}
\newcommand{\ec}{\mathfrak{c}}
\newcommand{\tec}{\tilde{\ec}}
\newcommand{\ef}{\mathfrak{f}}
\newcommand{\eh}{\mathfrak{h}}
\newcommand{\eo}{\mathfrak{o}}
\newcommand{\epp}{\mathfrak{p}}
\newcommand{\eX}{{\mathcal X}}
\newcommand{\oeX}{\overline{\eX}}
\newcommand{\oH}{\bar{H}}
\newcommand{\oX}{\overline{X}}
\newcommand{\ohne}{\setminus}
\newcommand{\silo}{\stackrel{\sim}{\longrightarrow}}
\newcommand{\tei}{\, | \,}
\newcommand{\ent}{\;\widehat{=}\;}
\newcommand{\hullet}{\raisebox{0.03cm}{$\scriptstyle \bullet$}}
\newcommand{\halb}{\frac{1}{2}}
\newcommand{\dis}{\displaystyle}
\newtheorem{theorem}{Theorem}

\newtheorem{prop}[theorem]{Proposition}
\newtheorem{remark}[theorem]{Remark}
\newtheorem{conjecture}[theorem]{Conjecture}
\newtheorem{dic}[theorem]{Dictionary}
\newtheorem{constr}[theorem]{Construction}
\begin{document}
\title{Analogies between analysis on foliated spaces and arithmetic geometry}
\author{Christopher Deninger}
\date{\ }
\maketitle
\section{Introduction} \label{s:1}
For the arithmetic study of varieties over finite fields powerful cohomological methods are available which in particular shed much light on the nature of the corresponding zeta functions. For algebraic schemes over $\spec \Z$ and in particular for the Riemann zeta function no cohomology theory has yet been developed that could serve similar purposes. For a long time it had even been a mystery how the formalism of such a theory could look like. This was clarified in \cite{D1}. However until now the conjectured cohomology has not been constructed. \\
There is a simple class of dynamical systems on foliated manifolds whose reduced leafwise cohomology has several of the expected structural properties of the desired cohomology for algebraic schemes. In this analogy, the case where the foliation has a dense leaf corresponds to the case where the algebraic scheme is flat over $\spec \Z$ e.g. to $\spec \Z$ itself. In this situation the foliation cohomology which in general is infinite dimensional is not of a topological but instead of a very analytic nature. This can also be seen from its description in terms of global differential forms which are harmonic along the leaves. An optimistic guess would be that for arithmetic schemes $\eX$ there exist foliated dynamical systems $X$ whose reduced leafwise cohomology gives the desired cohomology of $\eX$. If $\eX$ is an elliptic curve over a finite field this is indeed the case with $X$ a generalized solenoid, not a manifold, \cite{D3}. \\
We illustrate this philosophy by comparing the ``explicit formulas'' in analytic number theory to a transversal index theorem. We also give  a dynamical analogue of a recent conjecture of Lichtenbaum on special values of Hasse--Weil zeta functions by using the Cheeger-M\"uller theorem equating analytic and Reidemeister torsion. As a new example we point out an analogy between properties of Cramer's function made from zeroes of the Riemann zeta function and properties of the trace of a wave operator. Incidentally, in all cases analytic number theory suggests questions in the theory of partial differential operators on foliated manifolds. \\
Since the entire approach is not yet in a definitive state our style is deliberately a bit brief in places. For more details we refer to the references provided.\\
Hermann Weyl was very interested in number theory and one of the pioneers in the spectral theory of differential operators. In particular the theory of zeta-regularized determinants and their application to analytic torsion evolved from his work. In my opinion the analogies reviewed in the present note and in particular the analytic nature of foliation cohomology show that a deeper understanding of number theoretical zeta functions will ultimately require much more analysis of partial differential operators and of Weyl's work than presently conceived. 

I would like to thank the referee very much for thoughtful comments.

\section{Gradient flows}

We introduce zeta functions and give some arithmetic motivation for the class of dynamical systems considered in the sequel.

Consider the Riemann zeta function
\[
\zeta (s) = \prod_p (1 - p^{-s})^{-1} = \sum^{\infty}_{n=1} n^{-s} \quad \mbox{for} \; \RRe s > 1 \; .
\]
It has a holomorphic continuation to $\C \ohne \{ 1 \}$ with a simple pole at $s = 1$. The famous Riemann hypothesis asserts that all zeroes of $\zeta (s)$ in the critical strip $0 < \RRe s < 1$ should lie on the line $\RRe s = 1/2$. A natural generalization of $\zeta (s)$ to the context of arithmetic geometry is the Hasse--Weil zeta function $\zeta_{\eX} (s)$ of an algebraic scheme $\eX / \Z$
\[
\zeta_{\eX} (s) = \prod_{x \in |\eX|} (1 - Nx^{-s})^{-1} \quad \mbox{for} \; \RRe s > \dim \eX \; .
\]
Here $|\eX|$ is the set of closed points $x$ of $\eX$ and $Nx$ is the number of elements in the residue field $\kappa (x)$ of $x$. For the one-dimensional scheme $\eX = \spec \Z$ we recover $\zeta (s)$ and for $\eX = \spec \eo_K$ where $\eo_K$ is the ring of integers in a number field we obtain the Dedekind zeta function $\zeta_K (s)$. It is expected that $\zeta_{\eX} (s)$ has a meromorphic continuation to all of $\C$. This is known in many interesting cases but by no means in general.

If $\eX$ has characteristic $p$ then the datum of $\eX$ over $\F_p$ is equivalent to the pair $(\eX \otimes \oF_p , \varphi)$ where $\varphi = F \otimes \id_{\oF_p}$ is the $\oF_p$-linear base extension of the absolute Frobenius morphism $F$ of $\eX$. For example, the set of closed points $|\eX|$ of $\eX$ corresponds to the set of finite $\varphi$-orbits $\eo$ on the $\oF_p$-valued points $(\eX \otimes \oF_p) (\oF_p) = \eX (\oF_p)$ of $\eX \otimes \oF_p$. We have $\log N x = |\eo| \log p$ under this correspondence. Pairs $(\eX \otimes \oF_p , \varphi)$ are roughly analogous to pairs $(M, \varphi)$ where $M$ is a smooth manifold of dimension $2 \dim \eX$ and $\varphi$ is a diffeomorphism of $M$. A better analogy would be with K\"ahler manifolds and covering self maps of degree greater than one. See remark \ref{t33} below.

The $\Z$-action on $M$ via the powers of $\varphi$ can be suspended as follows to an $\R$-action on a new space. Consider the quotient
\begin{equation} \label{eq:31}
X = M \times_{p^{\Z}} \R^*_+
\end{equation}
where the subgroup $p^{\Z}$ of $\R^*_+$ acts on $M \times \R^*_+$ as follows:
\[
p^{\nu} \cdot (m , u) = (\varphi^{-\nu} (m) , p^{\nu} u) \quad \mbox{for} \; \nu \in \Z , m \in M , u \in \R^*_+ \; .
\]
The group $\R$ acts smoothly on $X$ by the formula
\[
\phi^t [m,u] = [m , e^t u] \; .
\]
Here $\phi^t$ is the diffeomorphism of $X$ corresponding to $t \in \R$. Note that the closed orbits $\gamma$ of the $\R$-action on $X$ are in bijection with the finite $\varphi$-orbits $\eo$ on $M$ in such a way that the length $l (\gamma)$ of $\gamma$ satisfies the relation $l (\gamma) = |\eo| \log p$. 

Thus in an analogy of $\eX / \F_p$ with $X$, the closed points $x$ of $\eX$ correspond to the closed orbits of the $\R$-action $\phi$ on $X$ and $Nx = p^{|\eo|}$ corresponds to $e^{l (\gamma)} = p^{|\eo|}$. Moreover if $d = \dim \eX$ then $\dim X = 2d + 1$. 

The system \eqref{eq:31} has more structure. The fibres of the natural projection of $X$ to $\R^*_+ / p^{\Z}$ form a $1$-codimensional foliation $\Fh$ (in fact a fibration). The leaves (fibres) of $\Fh$ are the images of $M$ under the immersions for every $u \in \R^*_+$ sending $m$ to $[m,u]$. In particular, the leaves are transversal to the flow lines of the $\R$-action and $\phi^t$ maps leaves to leaves for every $t \in \R$.

Now the basic idea is this: If $\eX / \Z$ is flat, there is no Frobenius and hence no analogy with a discrete time dynamical system i.e.~an action of $\Z$. However one obtains a reasonable analogy with a continuous time dynamical system by the following correspondence


\begin{dic}, {\bf part 1} \label{t32_1}
 
\em 
\begin{tabular}{p{7cm}|p{7cm}}
$\eX$ $d$-dimensional regular algebraic scheme over $\spec \Z$ & triple $(X , \Fh , \phi^t)$, where $X$ is a $2d+1$-dimensional smooth manifold with a smooth $\R$-action $\phi : \R \times X \to X$ and a one-codimensional foliation $\Fh$. The leaves of $\Fh$ should be transversal to the $\R$-orbits and every diffeomorphism $\phi^t$ should map leaves to leaves. \\ \hline
closed point $x$ of $\eX$ & closed orbit $\gamma$ of the $\R$ action \\ \hline
 Norm $Nx$ of closed point $x$ & $N \gamma = \exp l (\gamma)$ for closed orbit $\gamma$ \\ \hline
Hasse--Weil zeta function \newline
$\zeta_{\eX} (s) = \prod_{x \in |\eX|} (1 - Nx^{-s})^{-1}$ & Ruelle zeta function (\ref{eq:44}) \newline $\zeta_X (s) = \prod_{\gamma} (1 - N \gamma^{-s})^{\pm 1}$ \newline
(if the product makes sense) \\ \hline
$\eX \to \spec \F_p$ & triples $(X , \Fh , \phi^t)$ where $\Fh$ is given by the fibres of an $\R$-equivariant fibration $X \to \R^*_+ / p^{\Z}$
\end{tabular}
\end{dic}

\begin{remark} \label{t33}
\em 
A more accurate analogy can be motivated as follows: If $(M , \varphi)$ is a pair consisting of a manifold with a self covering $\varphi$ of degree $\deg \varphi \ge 2$ one can form the generalized solenoid $\hat{M} = \varprojlim (\ldots \xrightarrow{\varphi} M \xrightarrow{\varphi} M \xrightarrow{} \ldots)$ , c.f. \cite{CC} 11.3. On $\hat{M}$ the map $\varphi$ becomes the shift isomorphism and as before we may consider a suspension $X = \hat{M} \times_{p^{\Z}} \R^*_+$. This example suggests that more precisely schemes $\eX / \Z$ should correspond to triples $(X , \Fh, \phi^t)$ where $X$ is a $2d + 1$ dimensional generalized solenoid which is also a foliated space with a foliation $\Fh$ by K\"ahler manifolds of complex dimension $d$. It is possible to do analysis on such generalized spaces c.f. \cite{CC}. Analogies with arithmetic are studied in more detail in \cite{D2} and \cite{Le}. 
\end{remark}

\begin{constr} \label{t34}
\em
Let us consider a triple $(X , \Fh , \phi^t)$ as in the dictionary above, let $Y_{\phi}$ be the vector field generated by the flow $\phi^t$ and let $T \Fh$ be the tangent bundle to the foliation. Let $\omega_{\phi}$ be the one-form on $X$ defined by
\[
\omega_{\phi} \, |_{T\Fh} = 0 \quad \mbox{and} \quad \langle \omega_{\phi} , Y_{\phi} \rangle = 1 \; .
\]
One checks that $d\omega_{\phi} = 0$ and that $\omega_{\phi}$ is $\phi^t$-invariant i.e. $\phi^{t*} \omega_{\phi} = \omega_{\phi}$ for all $t \in \R$. We may view the cohomology class $\psi = [\omega_{\phi}]$ in $H^1 (X ,\R)$ defined by $\omega_{\phi}$ as a homomorphism
\[
 \psi : \pi^{\abb}_1 (X) \longrightarrow \R \; .
\]
Its image $\Lambda \subset \R$ is called the group of periods of $(X , \Fh , \phi^t)$.
\end{constr}

It is known that $\Fh$ is a fibration if and only if $\rank\, \Lambda = 1$. In this case there is an $\R$-equivariant fibration $X \longrightarrow \R / \Lambda$ whose fibres are the leaves of $\Fh$. Incidentally, a good reference for the dynamical systems we are considering is \cite{Fa}.

Consider another motivation for our foliation setting: For a number field $K$ and any $f \in K^*$ we have the product formula where $\epp$ runs over the places of $K$
\[
\prod_{\epp} \| f \|_{\epp} = 1 \; .
\]
Here for the finite places i.e. the prime ideals $\epp$ we have 
\[
\| f \|_{\epp} = N \epp^{-\ord_{\epp} f} \; .
\]
Now look at triples $(X , \Fh , \phi^t)$ where the leaves of $\Fh$ are Riemann surfaces varying smoothly in the transversal direction. In particular we have $\dim X = 3 = 2d + 1$ where $d  = 1 = \dim (\spec \eo_K)$. Consider smooth functions $f$ on $X$ which are meromorphic on leaves and have their divisors supported on closed orbits of the flow. For compact $X$ it follows from a formula of Ghys that
\[
\prod_{\gamma} \; \| f \|_{\gamma} = 1 \; .
\]
In the product $\gamma$ runs over the closed orbits and 
\[
\| f \|_{\gamma} = (N\gamma)^{-\ord_{\gamma}f}
\]
where $N\gamma = e^{l(\gamma)}$ and $\ord_{\gamma} f = \ord_z (f \, |_{F_z})$. Here $z$ is any point on $\gamma$ and $F_z$ is the leaf through $z$.\\
So the foliation setting allows for a product formula where the $N \gamma$ are not all powers of the same number. If one wants an infinitely generated $\Lambda$, one must allow the flow to have fixed points ($\ent$ infinite places). A product formula in this more general setting is given in \cite{Ko}.

We end this motivational section listing these and more analogies some of which will be explained later.

\begin{dic}, {\bf part 2} \label{t32_2} 

\em
\begin{tabular}{p{7cm}|p{7cm}}
Number of residue characteristics $|\car (\eX)|$ of $\eX$ & rank of period group $\Lambda$ of $(X, \Fh , \phi^t)$ \\ \hline
Weil \'etale cohomology of $\eX$ & Sheaf cohomology of $X$ \\ \hline 
\vspace*{0.1mm} Arakelov compactification $\oeX = \eX \cup \eX_{\infty}$ where $\eX_{\infty} = \eX \otimes \R$ & \vspace*{0.1mm} Triples $(\oX , \Fh , \phi^t)$ where $\oX$ is a $2d + 1$-dimensional compactification of $X$ with a flow $\phi$ and a $1$-codimensional foliation. The flow maps leaves to leaves, however $\phi$ may have fixed points on $\oX$\\ \hline
$\eX (\C) / F_{\infty}$ , where $F_{\infty}$ denotes complex conjugation & set of fixed points of $\phi^t$. Note that the leaf of $\Fh$ containing a fixed point is $\phi$-invariant. \\ \hline
%
$\spec \kappa (x) \hookrightarrow \eX$ for $x \in |\eX|$ & Embedded circle i.e. knot $\R / l (\gamma) \Z \hookrightarrow X$ corresponding to a periodic orbit $\gamma$ (map $t + l (\gamma) \Z$ to $\phi^t (x)$ for a chosen point $x$ of $\gamma$). \\ \hline
Product formula for number fields $\prod\limits_{\epp} \|f\|_{\epp} = 1$ & Kopei's product formula \cite{Ko}\\ \hline
Explicit formulas of analytic number theory & transversal index theorem for $\R$-action on $X$ and Laplacian along the leaves of $\Fh$, c.f. \cite{D3}  \\ \hline
$- \log |d_{K / \Q}|$ & Connes' Euler characteristic $\chi_{\Fh} (X , \mu)$ for Riemann surface laminations with respect to transversal measure $\mu$ defined by $\phi^t$, c.f. \cite{D2} section 4 \\ \hline
\end{tabular}
\end{dic}

\section{Explicit formulas and transversal index theory}
A simple version of the explicit formulas in number theory asserts the following equality of distributions in $\Dh' (\R^{> 0})$:
\begin{equation} \label{eq:1n}
1 - \sum_{\rho} e^{t\rho} + e^t = \sum_p \log p \sum_{k\ge 1} \delta_{k \log p} + (1 - e^{-2t})^{-1} \; .
\end{equation}
Here $\rho$ runs over the zeroes of $\zeta (s)$ in $0 < \RRe s < 1$, the Dirac distribution in $a \in \R$ is denoted by $\delta_a$ and the functions $e^{t \rho}$ are viewed as distributions. Note that in the space of distributions the sum $\sum_{\rho} e^{t\rho}$ converges as one sees by partial integration since $\sum_{\rho} \rho^{-2}$ converges.

We want to compare this formula with a transversal index theorem in geometry.

With our dictionary in mind consider triples $(X , \Fh , \phi^t)$ where $X$ is a smooth compact manifold of odd dimension $2d + 1$, equipped with a one-codimensional foliation $\Fh$ and $\phi^t$ is a flow mapping leaves of $\Fh$ to leaves. Moreover we assume that the flow lines meet the leaves transversally in every point so that $\phi$ has no fixed points. Consider the (reduced) foliation cohomology:
\[
\oH^{\hullet}_{\Fh} (X) := \Ker d_{\Fh} / \overline{\Imm d_{\Fh}} \; . 
\]
Here $(\Ah^{\hullet}_{\Fh} (X) , d_{\Fh})$ with $\Ah^i_{\Fh} (X) = \Gamma (X , \Lambda^i T^* \Fh)$ is the ``de Rham complex along the leaves'', (differentials only in the leaf direction). Moreover $\overline{\Imm d_{\Fh}}$ is the closure in the smooth topology of $\Ah^{\hullet}_{\Fh} (X)$, c.f. \cite{AK1}. 

The groups $\oH^{\hullet}_{\Fh} (X)$ have a smooth linear $\R$-action $\phi^*$ induced by the flow $\phi^t$. The infinitesimal generator $\theta = \lim_{t \to 0} \frac{1}{t} (\phi^{t*} - \id)$ exists on $\oH^{\hullet}_{\Fh} (X)$. It plays a similar role as the Frobenius morphism on \'etale or crystalline cohomology.

In general, the cohomologies $\oH^{\hullet}_{\Fh} (X)$ will be infinite dimensional Fr\'echet spaces. They are related to harmonic forms. Assume for simplicity that $X$ and $\Fh$ are oriented in a compatible way and choose a metric $g_{\Fh}$ on $T\Fh$. This gives a Hodge scalar product on $\Ah^{\hullet}_{\Fh} (X)$. We define the leafwise Laplace operator by
\[
\Delta_{\Fh} = d_{\Fh} d^*_{\Fh} + d^*_{\Fh} d_{\Fh} \quad \mbox{on} \; \Ah^{\hullet}_{\Fh} (X) \; .
\]

Then we have by a deep theorem by \'Alvarez L\'opez and Kordyukov \cite{AK1}:
\begin{equation}
  \label{eq:42}
  \oH^{\hullet}_{\Fh} (X) = \Ker \Delta_{\Fh} \; .
\end{equation}
Note that $\Delta_{\Fh}$ is not elliptic but only elliptic along the leaves of $\Fh$. Hence the standard regularity theory of elliptic operators does not suffice for \eqref{eq:42}. 

The isomorphism (\ref{eq:42}) is a consequence of the Hodge decomposition proved in \cite{AK1}
\begin{equation}
  \label{eq:43}
  \Ah^{\hullet}_{\Fh} (X) = \Ker \Delta_{\Fh} \oplus \overline{\Imm d_{\Fh} d^*_{\Fh}} \oplus \overline{\Imm d^*_{\Fh} d_{\Fh}} \; . 
\end{equation}
Consider $\Delta_{\Fh}$ as an unbounded operator on the space $\Ah^{\hullet}_{\Fh , (2)} (X)$ of leafwise forms which are $L^2$ on $X$. Then $\Delta_{\Fh}$ is symmetric and its closure $\overline{\Delta}_{\Fh}$ is selfadjoint. The orthogonal projection $P$ of $\Ah^{\hullet}_{\Fh, (2)} (X)$ to $\Ker \overline{\Delta}_{\Fh}$ restricts to the projection $P$ of $\Ah^{\hullet}_{\Fh} (X)$ to $\Ker \Delta_{\Fh}$ in \eqref{eq:43}. There is also an $L^2$-version $\oH^{\hullet}_{\Fh, (2)} (X)$ of the reduced leafwise cohomology and it follows from a result of von Neumann that we have
\[
\oH^{\hullet}_{\Fh , (2)} (X) = \Ker \overline{\Delta}_{\Fh} \; .
\]
In order to get a formula analogous to \eqref{eq:1n} we have to take $\dim X = 3$. Then $\Fh$ is a foliation by surfaces so that $\oH^i_{\Fh} (X) = 0$ for $i > 2$. The following transversal index theorem is contained in \cite{AK2} and \cite{DS}:

\begin{theorem}
Let $(X , \Fh , \phi^t)$ be as above with $\dim X = 3$ such that $\Fh$ has a dense leaf. Assume that there is a metric $g_{\Fh}$ on $T\Fh$ such that $\phi^t$ acts with conformal factor $e^{\alpha t}$. (This condition is very strong and can be replaced by weaker ones.) Then the spectrum of $\theta$ on $\oH^1_{\Fh, (2)} (X)$ consists of eigenvalues. If all periodic orbits of $\phi^t$ are non-degenerate we have an equality of distributions in $\Dh' (\R^{>0})$ for suitable (known) signs:
\begin{equation} \label{eq:2n}
1 - \sum_{\rho} e^{t\rho} + e^{\alpha t} = \sum_{\gamma} l (\gamma) \sum_{k \ge 1} \pm \delta_{kl (\gamma)} \; .
\end{equation}
In the sum $\rho$ runs over the spectrum of $\theta$ on $\oH^1_{\Fh, (2)} (X)$ and all $\rho$ satisfy $\RRe \rho = \frac{\alpha}{2}$.
\end{theorem}

Comparing \eqref{eq:1n} and \eqref{eq:2n} we see again that the primes $p$ and the closed orbits $\gamma$ should correspond with $\log p \ent l (\gamma)$. Since we assumed that the flow $\phi^t$ has no fixed points there is no contribution in \eqref{eq:2n} corresponding to the term $(1 - e^{-2t})^{-1}$ in \eqref{eq:1n}. Ongoing work of \'Alvarez L\'opez and Korduykov leads to optimism that there is an extension of the formula in the theorem to the case where $\phi^t$ may have fixed points and where the right hand side looks similar to \eqref{eq:1n}. However reduced leafwise $L^2$-cohomology may have to be replaced by ``adiabatic cohomology''.

The conditions in the theorem actually force $\alpha = 0$ so that $\phi^t$ is isometric. This is related to remark 2. In fact in \cite{D3} using elliptic curves over finite fields we constructed examples of systems $(X , \Fh , \phi^t)$ where $X$ is a solenoidal space and we have $\alpha = 1$. 

For simplicity we have only considered the explicit formula on $\R^{> 0}$. The comparison works also on all of $\R$ and among other features one gets a beautiful analogy between Connes' Euler characteristic $\chi_{\Fh} (X , \mu)$ and $- \log |d_{K / \Q}|$ c.f. \cite{D2} section 4.

\section{A conjecture of Lichtenbaum and a dynamical analogue}

Consider a regular scheme $\eX$ which is separated and of finite type over $\spec \Z$ and assume that $\zeta_{\eX} (s)$ has an analytic continuation to $s = 0$. 

Lichtenbaum conjectures the existence of a certain ``Weil-\'etale'' cohomology theory with and without compact supports, $H^i_c (\eX , A)$ and $H^i (\eX , A)$ for topological rings $A$. See \cite{Li1}, \cite{Li2}. It should be related to the zeta function of $\eX$ as follows.

\begin{conjecture}[Lichtenbaum] \label{t21}
Let $\eX / \Z$ be as above. Then the groups $H^i_c (\eX , \Z)$ are finitely generated and vanish for $i > 2 \dim \eX + 1$. Giving $\R$ the usual topology we have
\[
H^i_c (\eX , \Z) \otimes_{\Z} \R = H^i_c (\eX, \R) \; .
\]
Moreover, there is a canonical element $\psi$ in $H^1 (\eX,  \R)$ which is functorial in $\eX$ and such that we have:\\
{\bf a} The complex
\[
\ldots \xrightarrow{D} H^i_c (\eX , \R) \xrightarrow{D} H^{i+1}_c (\eX , \R) \xrightarrow{} \ldots
\]
where $Dh = \psi \cup h$, is acyclic. Note that $D^2 = 0$ because $\deg \psi = 1$.\\[0.2cm]
{\bf b} $\ord_{s=0} \zeta_{\eX} (s) = \sum_i (-1)^i i \, \rk \, H^i_c (\eX , \Z)$. \\[0.2cm]
{\bf c} For the leading coefficient $\zeta^*_{\eX} (0)$ of $\zeta_{\eX} (s)$ in the Taylor development at $s = 0$ we have the formula:
\[
\zeta^*_{\eX} (0) = \pm \prod_i |H^i_c (\eX, \Z)_{\tors} |^{(-1)^i} / \ddet (H^{\hullet}_c (\eX , \R) , D , \ef^{\hullet}) \; .
\]
Here, $\ef^i$ is a basis of $H^i_c (\eX , \Z) / \tors$.
\end{conjecture}

{\bf Explanation} For an acyclic complex of finite dimensional $\R$-vector spaces
\begin{equation} \label{eq:22}
0 \xrightarrow{} V^0 \xrightarrow{D} V^1 \xrightarrow{D} \ldots \xrightarrow{D} V^r \xrightarrow{} 0
\end{equation}
and bases $\eb^i$ of $V^i$ a determinant $\det (V^{\hullet} , D , \eb^{\hullet})$ in $\R^*_+$ is defined as follows: For bases $\ea = (w_1 , \ldots , w_n)$ and $\eb = (v_1 , \ldots , v_n)$ of a finite dimensional vector space $V$ set $[\eb / \ea] = \det M$ where $v_i = \sum_j m_{ij} w_j$ and $M = (m_{ij})$. Thus we have $[\ec / \ea] = [\ec / \eb] [\eb / \ea]$.

For all $i$ choose bases $\ec^i$ of $D (V^{i-1})$ in $V^i$ and a linearly independent set $\tilde{\ec}^{i-1}$ of vectors in $V^{i-1}$ with $D (\tilde{\ec}^{i-1}) = \ec^i$. Then $(\ec^i , \tec^i)$ is a basis of $V^i$ since \eqref{eq:22} is acyclic and one defines:
\begin{equation} \label{eq:24}
\det (V^{\hullet} , D , \eb^{\hullet}) = \prod_i |[\eb^i / (\ec^i , \tec^i)]|^{(-1)^i} \; .
\end{equation}
Note that $\det (V^{\hullet} , D , \eb^{\hullet})$ is unchanged if we replace the bases $\eb^i$ with bases $\ea^i$ such that $|[\eb^i / \ea^i]| = 1$ for all $i$. Thus the $\eb^i$ could be replaced by unimodularly or orthogonally equivalent bases. In particular the determinant
\[
\det (H^{\hullet}_c (\eX , \R) , D , \ef^{\hullet})
\]
does not depend on the choice of bases $\ef^i$ of $H^i_c (\eX , \Z) / \tors$. 

The following formula is obvious:
\begin{prop} \label{t25}
Let $\ea^i$ and $\eb^i$ be bases of the $V^i$ in \eqref{eq:22}. Then we have:
\[
\det (V^{\hullet} , D , \eb^{\hullet}) = \det (V^{\hullet} , D , \ea^{\hullet}) \prod_i |[\eb^i / \ea^i]|^{(-1)^i} \; .
\]
\end{prop}

For smooth projective varieties over finite fields using the Weil-\'etale topology Lichtenbaum has proved his conjecture in \cite{Li1}. See also \cite{Ge} for generalizations to the singular case. The formalism also works nicely in the study of $\zeta_{\eX} (s)$ at $s = 1/2$, c.f. \cite{Ra}. If $\eX$ is the spectrum of a number field, Lichtenbaum gave a definition of ``Weil-\'etale'' cohomology groups using the Weil group of the number field. Using the formula
\[
\zeta^*_K (0) = -\frac{hR}{w}
\]
he was able to verify his conjecture except for the vanishing of cohomology in degrees greater than three, \cite{Li2}. In fact, his cohomology does not vanish in higher degrees as was recently shown by Flach and Geisser so that some modification will be necessary. 

For a dynamical analogue let us look at a triple $(X , \Fh , \phi^t)$ as in section 2 with $X$ a closed manifold of dimension $2d + 1$ and $\phi^t$ everywhere transversal to $\Fh$. We then have a decomposition $TX = T \Fh \oplus T_0 X$ where $T_0 X$ is the rank one bundle of tangents to the flow lines.

In this situation the role of Lichtenbaum's Weil-\'etale cohomology is played by the ordinary singular cohomology with $\Z$ or $\R$-coefficients of $X$. Note that because $X$ is compact we do not have to worry about compact supports. From the arithmetic point of view we are dealing with a very simple analogue only!

Lichtenbaum's complex is replaced by
\[
(H^{\hullet} (X, \R) , D) \quad \mbox{where}\; Dh = \psi \cup h \quad \mbox{and} \; \psi = [\omega_{\phi}] \; .
\]
Here $\omega_{\phi}$ is the $1$-form introduced in construction 3 above.

Now assume that the closed orbits of the flow are non-degenerate. Then at least formally we have:
\begin{equation}
  \label{eq:44}
  \begin{array}{rcl}
    \zeta_R (s) & := & \dis \prod_{\gamma} (1 - e^{-sl (\gamma)})^{-\varepsilon_{\gamma}} \\
 & = & \dis \prod_i \ddet_{\infty} (s \cdot \id - \theta \tei \oH^i_{\Fh} (X))^{(-1)^{i+1}} \; .
  \end{array}
\end{equation}
Here $\gamma$ runs over the closed orbits, $\varepsilon_{\gamma} = \sgn \det (1 - T_x \phi^{l (\gamma)} \tei T_x \Fh)$ for any $x \in \gamma$, the Euler product converges in some right half plane and $\det_{\infty}$ is the zeta-regularized determinant. The functions $\ddet_{\infty} (s \cdot \id - \theta \tei \oH^i_{\Fh} (X))$ should be entire.

If the closed orbits of $\phi$ are degenerate one can define a Ruelle zeta function via Fuller indices \cite{Fu} and relation \eqref{eq:44} should still hold. In the present note we assume for simplicity that the action of the flow $\phi$ on $T\Fh$ is isometric with respect to $g_{\Fh}$ and we do not insist on the condition that the closed orbits should be non-degenerate. Then $\theta$ has pure eigenvalue spectrum with finite multiplicities on $\oH^{\hullet}_{\Fh} (X) = \Ker \Delta^{\hullet}_{\Fh}$ by \cite{DS}, Theorem 2.6. We define $\zeta_R (s)$ by the formula
\[
\zeta_R (s) = \prod_i \ddet_{\infty} (s \cdot \id - \theta \tei \oH^i_{\Fh} (X))^{(-1)^{i+1}} \; ,
\]
if the individual regularized determinants exist and define entire functions. The following result is proved in \cite{D4}. 



\begin{theorem}\label{t41}
Consider a triple $(X, \Fh, \phi^t)$ as above with $\dim X = 2d + 1$ and compatible orientations of $X$ and $\Fh$ such that the flow acts isometrically with respect to a metric $g_{\Fh}$ on $T\Fh$. For {\bf b} and {\bf c} assume that the above zeta-regularized determinants exist. Then the following assertions hold:\\
{\bf a} The complex
\[
\ldots \xrightarrow{D} H^i (X,\R) \xrightarrow{D} H^{i+1} (X , \R) \longrightarrow \ldots
\]
where $Dh = \psi \cup h$ with $\psi = [\omega_{\phi}]$ is acyclic.\\[0.2cm]
{\bf b} \quad $\ord_{s=0} \zeta_R (s) = \sum_i (-1)^i i \, \rk H^i (X, \Z)$.\\[0.2cm]
{\bf c} For the leading coefficient in the Taylor development at $s = 0$ we have the formula:
\[
\zeta^*_R (0) = \prod_i \, |H^i (X, \Z)_{\tors} |^{(-1)^i} / \det (H^{\hullet} (X,\R) , D, \ef^{\hullet}) \; .
\]
Here, $\ef^i$ is a basis of $H^i (X, \Z) / \tors$.
\end{theorem}

{\bf Idea of a proof of c} 
We define a metric $g$ on $X$ by $g = g_{\Fh} + g_0$ (orthogonal sum) on $TX = T\Fh \oplus T_0 X$ with $g_0$ defined by $\| Y_{\phi , x} \| = 1$ for all $x \in X$. Using techniques from the heat equation proof of the index theorem and assuming the existence of zeta-regularized determinants one can prove the following identity:
\begin{equation}
  \label{eq:412}
  \zeta^*_R (0) = T (X,g)^{-1} \; .
\end{equation}

Here $T(X,g)$ is the analytic torsion introduced by Ray and Singer using the spectral zeta functions of the Laplace operators $\Delta^i$ on $i$-forms on $X$
\begin{equation} \label{eq:413}
T (X,g) = \exp \sum_i (-1)^i \frac{i}{2} \zeta'_{\Delta^i} (0) \; .
\end{equation}
 
Equation \eqref{eq:412} is the special case of Fried's conjecture in \cite{Fr3}, see also \cite{Fr1}, \cite{Fr2} where the flow has an {\it integrable} complementary distribution. 
Next, we use the famous Cheeger--M\"uller theorem:
\begin{equation}
  \label{eq:416}
  T (X,g) = \tau (X,g) \qquad \mbox{c.f. \cite{Ch}, \cite{M}}
\end{equation}
where $\tau (X,g) = \tau (X,\eh_{\hullet})$ is the Reidemeister torsion with respect to the volume forms on homology given by the following bases $\eh_i$ of $H_i (X, \R)$. Choose  orthonormal bases $\eh^i$ of $\Ker \Delta^i$ with respect to the Hodge scalar product and view them as bases of $H^i (X, \R)$ via $\Ker \Delta^i \silo H^i (X, \R)$. Let $\eh_i$ be the dual base to $\eh^i$. 

If we are given two bases $\ea$ and $\eb$ of a real vector space we have the formula
\[
\tau (X,\eb_{\hullet}) = \tau (X , \ea_{\hullet}) \prod_i |[\eb_i / \ea_i]|^{(-1)^i}
\]
for any choices of bases $\ea_i , \eb_i$ of $H_i (X, \R)$. Let $\ef_{\hullet} , \ef^{\hullet}$ be dual bases of $H_{\hullet} (X, \Z) / \tors$ resp. $H^{\hullet} (X, \Z) / \tors$. Then we have in particular:
\begin{equation}
  \label{eq:417}
  \tau (X,\eh_{\hullet}) = \tau (X, \ef_{\hullet}) \prod_i |[\eh_i / \ef_i]|^{(-1)^i} \; .
\end{equation}
Now it follows from the definition of Reidemeister torsion that we have:
\[
\tau (X, \ef_{\hullet}) = \prod_i |H_i (X, \Z)_{\tors}|^{(-1)^i} \; .
\]
The cap-isomorphism:
\[
\_ \cap [X] : H^j (X, \Z) \silo H_{2d +1 -j} (X, \Z)
\]
now implies that
\begin{equation}
  \label{eq:418}
  \tau (X, \ef_{\hullet}) = \prod_i |H^i (X, \Z)_{\tors}|^{(-1)^{i+1}} \; .
\end{equation}
Next let us look at the acyclic complex
\[
(H^{\hullet} (X, \R) , D = \psi \cup \, \_ ) \cong (\Ker \Delta^{\hullet} , \omega_{\phi} \wedge \, \_ ) \; .
\]
Using the canonical decomposition into bidegrees
\[
\Ker \Delta^n = \omega_{\phi} \wedge (\Ker \Delta^{n-1}_{\Fh})^{\theta=0} \oplus (\Ker \Delta^n_{\Fh})^{\theta = 0}
\]
we see that it is isometrically isomorphic to 
\begin{equation}
  \label{eq:421}
  (M^{\hullet -1} \oplus M^{\hullet} , D)
\end{equation}
where $M^i = (\Ker \Delta^i_{\Fh})^{\theta = 0}$ and $D(m' , m) = (m,0)$. We may choose the orthonormal basis $\eh^i$ above to be of the form $\eh^i = (\omega_{\phi} \wedge \tilde{\eh}^{i-1}  , \tilde{\eh}^i)$ where $\tilde{\eh}^i$ is an orthonormal basis of $(\Ker \Delta^i_{\Fh})^{\theta= 0}$. For this basis, i.e. $(\tilde{\eh}^{i-1} , \tilde{\eh}^i)$ in the version (\ref{eq:421}) it is trivial from the definition (\ref{eq:24}) that
\[
|\det (H^{\hullet} (X, \R) , D , \eh^{\hullet})| = 1 \; .
\]
Using proposition \ref{t25} for $\eb^{\hullet} = \eh^{\hullet}$ and $\ea^{\hullet} = \ef^{\hullet}$ we find
\begin{equation}
  \label{eq:422}
  1 = |\det (H^{\hullet} (X, \R) , D , \ef^{\hullet})| \prod_i |[\eh^i / \ef^i]|^{(-1)^i} \; .
\end{equation}
Hence we get
\begin{eqnarray*}
  \lefteqn{\prod_i |H^i (X, \Z)_{\tors}|^{(-1)^i} / \det (H^{\hullet} (X, \R) , D, \ef^{\hullet})}\\
& \overset{(\ref{eq:422})}{=} & \prod_i |H^i (X, \Z)_{\tors} |^{(-1)^i} \prod_i |[\eh^i / \ef^i|]^{(-1)^i} \\
& \overset{(\ref{eq:418})}{=} & \tau (X, \ef_{\hullet})^{-1} \prod_i |[ \eh_i / \ef_i |]^{(-1)^{i+1}} \overset{(\ref{eq:417})}{=} \tau (X, \eh_{\hullet})^{-1} \\
& \overset{(\ref{eq:416})}{=} & T (X, g)^{-1} \overset{(\ref{eq:412})}{=} \zeta^*_R (0) \; .
\end{eqnarray*}
\hspace*{\fill} $\Box$

{\bf The simplest example} For $q > 1$ let $X$ be the circle $\R / (\log q) \Z$ foliated by points and with $\R$ acting by translation. The Ruelle zeta function is given by
\[
\zeta_R (s) = (1 - e^{-s \log q})^{-1} = (1 - q^{-s})^{-1} \; .
\]
The operator $\theta = d / dx$ on $\oH^0_{\Fh} (X) = C^{\infty} (X)$ defines an unbounded operator on $\oH^0_{\Fh, (2)} (X) = L^2 (X)$ with pure eigenvalue spectrum $2 \pi i \nu / \log q$ for $\nu \in \Z$. It follows from a formula of Lerch e.g. \cite{D1} \S\,2 that we have
\[
\ddet_{\infty} (s \cdot \id - \theta \tei \oH^0_{\Fh} (X)) = 1 - q^{-s}
\]
and hence formula \eqref{eq:44} holds in this case.

Next, we have $\omega_{\phi} = dt$ and hence $\psi = [dt]$ is a generator of $H^1 (X , \R)$. The complex
\[
\ldots \xrightarrow{D} H^i (X, \R) \xrightarrow{D} H^{i+1} (X, \R) \to \ldots
\]
is therefore acyclic. Moreover we have
\[
\ord_{s = 0} \zeta_R (s) = -1 = \sum_i (-1)^i i \, \rk H^i (X , \Z) \; .
\]
The leading coefficient $\zeta^*_R (0)$ is given by $\zeta^*_R (0) = (\log q)^{-1}$. We have $H^i (X , \Z)_{\tors} = 0$ for all $i$. As $\psi / \log q$ is a basis of $H^1 (X , \Z)$ the formula $D (1) = \psi = (\log q) (\psi / \log q)$ for $1 \in H^0 (X , \R)$ shows that we have
\[
\det (H^{\hullet} (X , \R) , D , \ef^{\hullet}) = \log q \; .
\]
This illustrates part {\bf c} of theorem \ref{t41}.

\section{Cram\'ers function and the transversal wave equation}
In this section we observe a new analogy and mention some further directions of research in analysis suggested by this analogy. Let us write the non-trivial zeroes of $\zeta (s)$ as $\rho = \halb + \gamma$. In \cite{C} Cram\'er studied the function $W (z)$ which is defined for $\Imm z > 0$ by the absolutely and locally uniformly convergent series
\[
W (z) = \sum_{\Imm \gamma > 0} e^{\gamma z} \; .
\]
According to Cram\'er the function $W$ has a meromorphic continuation to $\C_- = \C \ohne \{ xi \tei x \le 0 \}$ with poles only for $z = m \log p$ with $m \in \Z , m \neq 0$ and $p$ a prime number. The poles are of first order. If we view the locally integrable function $e^{\gamma t}$ of $t$ as a distribution on $\R$ the series
\[
W_{\ddis} = \sum_{\Imm \gamma > 0} e^{\gamma t}
\]
converges in $\Dh' (\R)$. Using the mentioned results of Cram\'er on $W$ one sees that the singular support of $W_{\ddis}$ consists of $t = 0$ and the numbers $t = m \log p$ for $p$ a prime and $m \in \Z , m \neq 0$, c.f. \cite{DSch} section 1. One my ask if there is an analytic counterpart to this result in the spirit of the dictionary in section 2. 

Presently this is not the case although the work \cite{K} is related to the problem. According to our dictionary we consider a system $(X , \Fh , \phi^t)$ with a $3$-manifold $X$. Let us first assume that $X$ is compact and that $\phi^t$ is isometric and has no fixed points. Then we have $-\theta^2 = \Delta_0$ on $\Ah^{\hullet}_{\Fh} (X)$ where $\theta = \lim_{t \to 0} \frac{1}{t} (\phi^{t*} - \id)$ is the infinitesimal generator of the induced group of operators $\phi^{t*}$ on the Fr\'echet space $\Ah^{\hullet}_{\Fh} (X)$. Moreover $\Delta_0$ is the Laplacian along the flow lines with coefficients in $\Lambda^{\hullet} T^* \Fh$. Note that it is transversally elliptic with respect to $\Fh$. Since $\Delta_0 \, |_{\Ker \Delta_{\Fh}} = \Delta \, |_{\Ker \Delta_{\Fh}}$ by isometry, it follows from the corresponding result for the spectrum of the Laplacian that the spectrum $\{ \gamma \}$ of $\theta$ on $\oH^1_{\Fh , (2)} (X) = \Ker \overline{\Delta}^1_{\Fh}$ consists of eigenvalues and that all $\gamma$ are purely imaginary. (This proves a part of theorem 5.) For $\varphi$ in $\Dh (\R)$ the operator $\int_{\R} \varphi (t) e^{it |\theta|} \, dt$ on $\Ker \overline{\Delta}^1_{\Fh}$ is trace class and the map $\varphi \mapsto \tr \int_{\R} \varphi (t) e^{it |\theta|}$ defines a distribution denoted by $\tr (e^{it |\theta|} \tei \Ker \overline{\Delta}^1_{\Fh})$. Then we have:
\begin{eqnarray*}
2 \sum_{\Imm \gamma > 0} e^{t \gamma} & = & \tr (e^{it |\theta|} \tei \Ker \overline{\Delta}^1_{\Fh}) \\
 & = & \tr (e^{it \sqrt{\Delta}} \tei \Ker \overline{\Delta}^1_{\Fh}) \quad \mbox{in} \; \Dh' (\R) \; .
\end{eqnarray*}
Here we have assumed that zero is not in the spectrum of $\theta$ on $\Ker \overline{\Delta}^1_{\Fh}$, corresponding to the fact that $\zeta (1/2) \neq 0$. But this is not important. 

On functions, i.e. on $ C^{\infty} (X) = \Ah^0_{\Fh} (X)$ instead of $\Ker \Delta^1_{\Fh}$ the distributional trace $\tr (e^{it \sqrt{\Delta}})$ of the operator $e^{it \sqrt{\Delta}}$ has been extensively studied. By a basic result of Chazarain \cite{Cha} the singular support of $\tr (e^{it \sqrt{\Delta}})$ consists of $t = 0$ and the numbers $t = m l (\gamma)$ for $m \in \Z , m \neq 0$ and $l (\gamma)$ the lengths of the closed orbits of the geodesic flow. On the other hand by the analogy with Cram\'ers theorem we expect that the singularities of $\tr (e^{it \sqrt{\Delta}} \tei \Ker \overline{\Delta}^1_{\Fh})$ are contained in the analogous set where now $\gamma$ runs over the closed orbits of the flow $\phi^t$. 

Based on the work of H\"ormander the ``big'' singularity of $\tr (e^{it \sqrt{\Delta}})$ at $t = 0$ was analyzed in \cite{DG} proposition 2.1 by an asymptotic expansion of the Fourier transform of $\tr (e^{it \sqrt{\Delta}})$. An analogue of this expansion for the Fourier transform of $\tr (e^{it \sqrt{\Delta}} \tei \Ker \overline{\Delta}^1_{\Fh})$ would correspond well with asymptotics that can be obtained with some effort from Cram\'er's theory, except that in Cram\'er's theory there also appear logarithmic terms coming from the infinite place. These should also appear in the analysis of $\tr (e^{it \sqrt{\Delta}} \tei \Ker \overline{\Delta}^1_{\Fh})$ if one allows the flow to have fixed points. Finally one should drop the condition that the flow acts isometrically. Then one has to work with $\tr (e^{it |\theta|} \tei \oH^1_{\Fh , (2)} (X))$ instead of $\tr (e^{it \sqrt{\Delta}} \tei \Ker \overline{\Delta}^1_{\Fh})$. Also one can try to prove a Duistermaat--Guillemin trace formula for $\tr (e^{it |\theta|} \tei \oH^1_{\Fh , (2)} (X))$ and even if the flow has fixed points.

Mathematisches Institut\\
Einsteinstr. 62\\
48149 M\"unster, Germany\\
c.deninger@math.uni-muenster.de
\end{document}